\documentclass{amsart}

\theoremstyle{definition}

\theoremstyle{remark}

\reversemarginpar

\begin{document}

\theoremstyle{plain}
\newtheorem{thm}{Theorem}[section]
\newtheorem{prop}{Proposition}[section]
\newtheorem{lem}{Lemma}[section]
\newtheorem{cor}{Corollary}[section]
\theoremstyle{definition}
\newtheorem{exm}{Example}[section]
\newtheorem{rem}{Remark}[section]
\newtheorem{dfn}{Definition}[section]
\numberwithin{equation}{section}
\renewcommand{\setminus}{\smallsetminus}
\newcommand{\s}{\vspace{0.2cm}}

\newtheorem{Definition}{\bf Definition}[section]
\newtheorem{Thm}[Definition]{\bf Theorem} 
\newtheorem{Example}[Definition]{\bf Example} 
\newtheorem{Lem}[Definition]{\bf Lemma} 
\newtheorem{Note}[Definition]{\bf Note} 
\newtheorem{Cor}[Definition]{\bf Corollary} 
\newtheorem{Prop}[Definition]{\bf Proposition} 
\newtheorem{Problem}[Definition]{\bf Problem} 
\numberwithin{equation}{section}

\title[Logarithms, Rational and Trigonometric]{The 
integrals in Gradshteyn and Ryzhik. \\
Part 9: Combinations of logarithms, rational and trigonometric functions.}

\subjclass[2000]{Primary 33}



\author[T. Amdeberhan]{Tewodros Amdeberhan}
\address{Department of Mathematics,
Tulane University, New Orleans, LA 70118}
\email{tamdeber@tulane.edu}

\author[V. Moll]{Victor H. Moll}
\address{Department of Mathematics,
Tulane University, New Orleans, LA 70118}
\email{vhm@math.tulane.edu}

\author[J. Rosenberg]{Jason Rosenberg}
\address{Department of Mathematics,
Tulane University, New Orleans, LA 70118}
\email{jrosenbe@tulane.edu}

\author[A. Straub]{Armin Straub}
\address{Department of Mathematics,
Tulane University, New Orleans, LA 70118}
\email{astraub@math.tulane.edu}

\author[P. Whitworth]{Pat Whitworth}
\address{Department of Mathematics,
Tulane University, New Orleans, LA 70118}
\email{pwhitwor@tulane.edu}

\begin{abstract}
The table of Gradshteyn and Ryzhik contains many integrals with integrands
of the form $R_{1}(x) \, \left( \ln R_{2}(x) \right)^{m}$, where 
$R_{1}$ and $R_{2}$ are
rational functions. In this paper we describe some examples where the 
logarithm appears to a single power, that is $m=1$, and the poles of 
$R_{1}$ are either real or purely imaginary.
\end{abstract}

\thanks{The authors  wish to the 
partial support of NSF-DMS 0409968 and NSF-CCLI 0633223.}

\maketitle

\newcommand{\nn}{\nonumber}
\newcommand{\ba}{\begin{eqnarray}}
\newcommand{\ea}{\end{eqnarray}}
\newcommand{\ift}{\int_{0}^{\infty}}
\newcommand{\ione}{\int_{0}^{1}}
\newcommand{\ifft}{\int_{- \infty}^{\infty}}
\newcommand{\no}{\noindent}
\newcommand{\realpart}{\mathop{\rm Re}\nolimits}
\newcommand{\imagpart}{\mathop{\rm Im}\nolimits}

\section{Introduction} \label{intro} 
\setcounter{equation}{0}

The table of integrals \cite{gr} contains many examples of the form 
\begin{equation}
\int_{a}^{b} R_{1}(x) \, \left( \ln R_{2}(x) \right)^{m} \, dx,
\label{gen-type}
\end{equation}
\noindent
where $R_{1}$ and $R_{2}$ are rational functions, 
$a, \, b \in \mathbb{R}^{+}$ and $m 
\in \mathbb{N}$. For example, $\mathbf{4.231.1}$ states that
\begin{equation}
\int_{0}^{1} \frac{\ln x \,dx }{1+x} = - \frac{\pi^{2}}{12}.
\end{equation}
\noindent
This result can be established by the elementary methods described here.

Other examples, such as  $\mathbf{4.233.1}$ 
\begin{equation}
\int_{0}^{1} \frac{\ln x \, dx }{1+x+x^{2}} = 
\frac{2}{9} \left( \frac{2 \pi^{2}}{3} - \psi'(\tfrac{1}{3}) \right),
\end{equation}
\noindent
and $\mathbf{4.261.8}$
\begin{equation}
\int_{0}^{1} \ln^{2}x \, \frac{1-x}{1-x^{6}} \, dx = 
\frac{8 \sqrt{3} \pi^{3} + 351 \zeta(3)}{486},
\end{equation}
\noindent
require more sophisticated special functions. Here 
$\psi$ is the {\em digamma function} defined by 
\begin{equation}
\psi(x) = \frac{\Gamma'(x)}{\Gamma(x)},
\end{equation}
\noindent
and $\zeta(s)$ is the classical 
{\em Riemann zeta function}. These results will be described in a future 
publication.  

The integrals discussed here can also be framed in the context of 
trigonometric functions. For example, the change of variables $x = \tan t$
yields the identity
\begin{equation}
\int_{0}^{1} \frac{\ln x \, dx}{1+x^{2}} = \int_{0}^{\pi/4} 
\ln \tan t \, dt = -G.
\label{trig1}
\end{equation}
\noindent
Here $G$ is the {\em Catalan's constant} defined by
\begin{equation}
G := \sum_{k=1}^{\infty} \frac{(-1)^{k}}{(2k+1)^{2}}. 
\end{equation}

In this paper we concentrate on integrals of the type (\ref{gen-type}) where
the logarithm appears to the first power and the poles of the rational
function are either real of purely imaginary.  The method of partial fractions
and scaling of the independent variable, show that such 
integrals are linear combinations of 
\begin{equation}
h_{n,1}(b) := \int_{0}^{b} \frac{\ln t \, dt}{(1+t)^{n}},
\end{equation}
\noindent
and 
\begin{equation}
h_{n,2}(b) := \int_{0}^{b} \frac{\ln t \, dt}{(1+t^{2})^{n}}. 
\end{equation}
\noindent
The function $h_{n,1}$ was evaluated in \cite{moll-gr2}, where is was 
denoted simply by $h$. We complete this 
evaluation in Section \ref{sec-singlepole}, by 
identifying a polynomial defined in \cite{moll-gr2}. The
closed-form of $h_{n,1}$ involves the {\em Stirling numbers} of the first 
kind. The evaluation of $h_{n,2}$ is discussed in Section  \ref{sec-imag}.
The value of $h_{n,2}$ involves the {\em tangent integral}
\begin{equation}
\text{Ti}_{2}(x) := \int_{0}^{x} \frac{\tan^{-1}t}{t}. 
\end{equation}

The case of integrals with more complicated pole structure will be 
described in a future publication.

\section{Some elementary examples} \label{sec-elementary} 
\setcounter{equation}{0}

We begin our discussion with an elementary example. Entry $\mathbf{4.291.1}$
states that
\begin{equation}
\int_{0}^{1} \frac{\ln(1+x)}{x} \,dx = \frac{\pi^{2}}{12}.
\label{elem1}
\end{equation}
\noindent
To establish this value we consider first a useful series. The result is expressed in terms
of the {\em Riemann zeta function}
\begin{equation}
\zeta(s) := \sum_{k=1}^{\infty} \frac{1}{k^{s}}.
\end{equation}

\begin{Lem}
Let $s >1$. Then 
\begin{equation}
\sum_{k=1}^{\infty} \frac{(-1)^{k}}{k^{s}} = -\frac{2^{s-1}-1}{2^{s-1}} \zeta(s)
\end{equation}
\noindent
and 
\begin{equation}
\sum_{k=1}^{\infty} \frac{1}{(2k-1)^{s}} = \frac{2^{s}-1}{2^{s}} \zeta(s)
\end{equation}
\end{Lem}
\begin{proof}
The second sum is 
\begin{equation}
\sum_{k=1}^{\infty} \frac{1}{(2k-1)^{s}} = \sum_{k=1}^{\infty} \frac{1}{k^{s}} -
\sum_{k=1}^{\infty} \frac{1}{(2k)^{s}} = (1-2^{-s}) \zeta(s). 
\end{equation}
\noindent
To evaluate the first sum, split it into even and odd values of the index $k$:
\begin{equation}
\sum_{k=1}^{\infty} \frac{(-1)^{k}}{k^{s}} = 
\sum_{k=1}^{\infty} \frac{1}{(2k)^{s}} - \sum_{k=1}^{\infty} \frac{1}{(2k-1)^{s}}
\end{equation}
\noindent
and use the value of the first sum. 
\end{proof}

To evaluate (\ref{elem1}) we employ the expansion
\begin{equation}
\ln(1+x) = \sum_{k=1}^{\infty} \frac{(-1)^{k-1}}{k}x^{k}
\end{equation}
\noindent
and integrate term by term we obtain
\begin{equation}
\int_{0}^{1} \frac{\ln(1+x)}{x} \,dx =  - \sum_{k=1}^{\infty} \frac{(-1)^{k}}{k^{2}}. 
\end{equation}
\noindent
The result now follows from the lemma and the classical value $\zeta(2) = \pi^{2}/6$. \\

A similar calculation yields $\mathbf{4.291.2}$:
\begin{equation}
\int_{0}^{1} \frac{\ln(1-x)}{x} dx = - \frac{\pi^{2}}{6}.
\end{equation}
\noindent
The change of variables $x = e^{-t}$ produces the evaluation of $\mathbf{4.223.1}$
\begin{equation}
\ift \ln(1 + e^{-t}) \, dt = \frac{\pi^{2}}{12} 
\end{equation}
\noindent
and $\mathbf{4.223.2}$:
\begin{equation}
\ift \ln(1 - e^{-t}) \, dt = -\frac{\pi^{2}}{6} 
\end{equation}

\section{More elementary examples} \label{sec-elem2} 
\setcounter{equation}{0}

In \cite{moll-gr2} we analyze the case in which the rational function $R_{1}$ 
has a single multiple pole and $R_{2}(x) = x$.  There are simple examples 
where the evaluation can be obtained directly. For instance, formula 
$\mathbf{4.231.1}$ states that
\begin{equation}
\int_{0}^{1} \frac{\ln x}{1+x} \, dx = - \frac{\pi^{2}}{12}.
\label{form888}
\end{equation}
\noindent
This can be established by simply expading the term $1/(1+x)$ as a geometric 
series and integrate term by term. The same is true for $\mathbf{4.231.2}$
\begin{equation}
\int_{0}^{1} \frac{\ln x}{1-x} \, dx = - \frac{\pi^{2}}{6}.
\label{form88}
\end{equation}
\noindent
The evaluation of $\mathbf{4.231.3}$
\begin{equation}
\int_{0}^{1} \frac{x \, \ln x}{1-x} \, dx = 1 - \frac{\pi^{2}}{6},
\end{equation}
\noindent
and $\mathbf{4.231.4}$ 
\begin{equation}
\int_{0}^{1} \frac{1+x}{1-x}  \, \ln x \, dx = 1 - \frac{\pi^{2}}{3},
\end{equation}
\noindent
follow directly from \ref{form88}. Similar elementary algebraic manipulations
produce $\mathbf{4.231.19}$
\begin{equation}
\int_{0}^{1} \frac{x \, \ln x}{1+x} \, dx = -1 + \frac{\pi^{2}}{2},
\end{equation}
\noindent
and $\mathbf{4.231.20}$
\begin{equation}
\int_{0}^{1} \frac{(1-x) \, \ln x}{1+x} \, dx = 1 - \frac{\pi^{2}}{6}.
\end{equation}

The same is true for $\mathbf{4.231.14}$
\begin{equation}
\int_{0}^{1} \frac{x \, \ln x }{1+x^{2}} \, dx = - \frac{\pi^{2}}{48},
\end{equation}
\noindent
and $\mathbf{4.231.15}$
\begin{equation}
\int_{0}^{1} \frac{x \, \ln x }{1-x^{2}} \, dx = - \frac{\pi^{2}}{24},
\end{equation}
\noindent
via the change of variables $t = x^{2}$. \\

The evaluation of $\mathbf{4.231.13}$:
\begin{equation}
\int_{0}^{1} \frac{\ln x \, dx}{1-x^{2}} = - \frac{\pi^{2}}{48}
\end{equation}
\noindent
is a direct consequence of the partial fraction decomposition 
\begin{equation}
\frac{1}{1-x^{2}} = \frac{1}{2} \frac{1}{1-x} + \frac{1}{2} \frac{1}{1+x}
\end{equation}
\noindent
and the identites (\ref{form888}) and (\ref{form88}). \\

It is often the  case that a simple change of variables reduces an integral 
to one that has previously been evaluated. For example, the change of variables 
$t = 1-x^{2}$ produces 
\begin{equation}
\int_{0}^{1} \frac{\ln(1-x^{2})}{x} \,  dx = \frac{1}{2} 
\int_{0}^{1} \frac{\ln t \,  dt}{1-t}.
\end{equation}
\noindent
The right-hand side is given in (\ref{form88}) and we obtain the
value of $\mathbf{4.295.11}$:
\begin{equation}
\int_{0}^{1} \frac{\ln(1-x^{2}) }{x} \, dx = - \frac{\pi^{2}}{12}.
\end{equation}

\section{A single multiple pole} \label{sec-singlepole} 
\setcounter{equation}{0}

The situation for a single multiple pole is more delicate. The pole 
pole may be placed at $x=-1$
by scaling. The main result established in \cite{moll-gr2} is: 

\begin{Thm}
Let $n \geq 2$ and $b >0$. Define 
\begin{equation}
h_{n,1}(b) = \int_{0}^{b} \frac{\ln t \, dt}{(1+t)^{n}}
\label{h1-def}
\end{equation}
\noindent
and introduce the function
\begin{equation}
q_{n}(b) = (1+b)^{n-1}h_{n,1}(b). 
\end{equation}
\noindent
Then 
\begin{equation}
q_{n}(b) = X_{n}(b) \ln b + Y_{n}(b) \ln(1+b) + Z_{n}(b),
\end{equation}
\noindent
where 
\begin{equation}
X_{n}(b) = \frac{(1+b)^{n-1}-1}{n-1}, \quad 
Y_{n}(b) = - \frac{(1+b)^{n-1}}{n-1}.
\end{equation}
\noindent
Finally, the function 
\begin{equation}
T_{n}(b) := - \frac{(n-1)! \, Z_{n}(b)}{b(1+b)},
\end{equation}
\noindent
satisfies $T_{2}(b) = 0$ and for $n \geq 1$ it satisfies the recurrence 
\begin{equation}
T_{n+2}(b) = n(1+b)T_{n+1}(b) + (n-1)! \left( \frac{(1+b)^{n}-1}{b} \right).
\label{recurT}
\end{equation}
\noindent
It follows that $T_{n}(b)$ is a polynomial in $b$ of degree $n-3$ with positive
integer coefficients.
\end{Thm}

\begin{Note}
The case $n=1$ is expressed in terms of the {\em polylogarithm function}
\begin{equation}
\text{PolyLog}[n,z] := \sum_{k=1}^{\infty} \frac{z^{k}}{k^{n}}.
\end{equation}
\noindent
Indeed, we have
\begin{equation}
\int_{0}^{b} \frac{\ln x}{1+x} \, dx = \ln b \, \ln(1+b) + \text{PolyLog}[2,-b].
\end{equation}
\end{Note}

\medskip

We now identify the polynomial $T_{n}$. The first few are given by
\begin{eqnarray}
T_{3}(b) & = & 1, \label{valuesT} \\
T_{4}(b) & = & 3b+4, \nonumber \\
T_{5}(b) & = & 11b^{2}+27b+18, \nonumber \\
T_{6}(b) & = & 50b^{3} + 176b^{2}+216b+96. \nonumber
\end{eqnarray}

For $n \geq 2$, define
\begin{equation}
A_{n}(b) = \frac{1}{(n-2)!}T_{n}(b).
\end{equation}
\noindent
Then (\ref{recurT}) becomes
\begin{equation}
A_{n+2}(b) = (1+b)A_{n+1}(b) + \frac{(1+b)^{n}-1}{bn},
\label{recurA}
\end{equation}
\noindent
with initial condition $A_{2}(b) = 0$. \\

The polynomial $A_{n}$ is written as
\begin{equation}
A_{n}(b) = \sum_{j=0}^{n-3} a_{n,j}b^{j}
\end{equation}
\noindent
The recursion (\ref{recurA}) produces 

\begin{Lem}
Let $n \geq 4$. Then the coefficients $a_{n,j}$ satisfy
\begin{eqnarray}
a_{n,0} & = & a_{n-1,0} + 1, \label{recuraa} \\
a_{n,j} & = & a_{n-1,j} + a_{n-1,j-1} + \frac{(n-3)!}{(j+1)! \, (n-3-j)!}, \text{ for } 
1 \leq j \leq n-4, \nonumber \\
a_{n,n-3} & = & a_{n-1,n-4} + \frac{1}{n-2}. \nonumber
\end{eqnarray}
\end{Lem}

The expressions  $a_{n,0} = n-2$ and $a_{n,n-3} = H_{n-2}$, the 
harmonic number, are easy to
determine from (\ref{recuraa}). We now find
closed-form expressions for the remaining 
coefficients. These involve the {\em Stirling numbers of the first kind}
$s(n,j)$ defined by the expansion
\begin{equation}
(x)_{n} = \sum_{j=1}^{n} s(n,j) x^{j},
\end{equation}
\noindent
where $(x)_{n} := x(x+1)(x+2) \cdots (x+n-1)$ is the Pochhammer symbol (also
called rising factorial). The numbers $s(n,1)$ are given by
\begin{equation}
s(n,1) = (-1)^{n-1} (n-1)!,
\end{equation}
\noindent
and the sequence $s(n,j)$ satisfies the recurrence 
\begin{equation}
\label{recu-stir}
s(n+1,j) = s(n,j-1) - n s(n,j).
\end{equation}

\begin{Thm}
The coefficients $a_{n,j}$ are given by
\begin{equation}
a_{n,j} = \frac{(-1)^{j}}{(j+1)!}  \binom{n-2}{j+1} \, s(j+2,2).
\label{defsn}
\end{equation}
\end{Thm}
\begin{proof}
Define 
\begin{equation}
b_{n,j} := (-1)^{j} a_{n,j} \times (j+1)! 
\binom{n-2}{j+1}^{-1}.
\end{equation}
\noindent
The recurrence (\ref{recuraa})  shows that $b_{n,j}$ is independent of $n$
and satisfies 
\begin{equation}
b_{n,j+1} = -j b_{n,j} +(-1)^{j+1}j!
\end{equation}
\noindent
and this is (\ref{recu-stir}) in the special case $j=2$. Formula (\ref{defsn})
has been established. 
\end{proof}

We now restate the value of $h_{n,1}$. The index $n$ is increased by $1$
and the identity 
$|s(n,k)| = (-1)^{n+k} s(n,k)$ is used
in order to make the formula look cleaner.

\begin{Cor}
\label{loga11}
Assume $ b > 0 $ and $n \in \mathbb{R}$. Then 
\begin{eqnarray}
\int_{0}^{b} \frac{\ln t \, dt}{(1+t)^{n+1}} & =  &  
\frac{1}{n} \left[ 1 - (1+b)^{-n} \right] \ln b - \frac{1}{n}  \ln(1+b)
\label{stir1}  \\
& - & \frac{1}{n (1+b)^{n-1}} \sum_{j=1}^{n-1} 
\frac{1}{j!} \binom{n-1}{j} |s(j+1,2)| b^{j}. \nonumber
\end{eqnarray}
\end{Cor}

\medskip

The special case $b=1$ provides the evaluation
\begin{equation}
\int_{0}^{1} \frac{\ln t \, dt}{(1+t)^{n+1}}  =
- \frac{\ln 2}{n}  
 -  \frac{1}{n 2^{n-1}} \sum_{j=1}^{n-1} 
\frac{1}{j!} \binom{n-1}{j} |s(j+1,2)|. \nonumber
\end{equation}

\medskip

Elementary changes of variables, starting with $t = \tan^{2} \varphi$, convert
(\ref{stir1}) into
\begin{eqnarray}
\int_{a}^{1} s^{n-1} \ln(1-s) \, ds & = & \frac{(1-a^{n})}{n^{2}} 
\left[ n \ln(1-a) -1 \right] \label{stir2} \\
& - & \frac{1}{n} \sum_{j=1}^{n-1} 
\frac{1}{j!}  \binom{n-1}{j} |s(j+1,2)| a^{n+1-j} (1-a)^{j}. \nonumber
\end{eqnarray}
\noindent
The special case $a=1/2$ produces
\begin{equation}
\int_{1/2}^{1} s^{n-1} \ln(1-s) \,ds = 
\frac{1}{n2^{n-1}} \sum_{j=1}^{n-1} \frac{\binom{n-1}{j} | s(j+1,2)| }{j!}
- \frac{(2^{n}-1)}{n^{2} 2^{n}} ( n \ln 2 + 1),
\nonumber
\end{equation}
\noindent
and $a=0$ gives
\begin{equation}
\int_{0}^{1} s^{n-1} \ln(1-s) \, ds = 
- \left( \frac{1}{n^{2}} + \frac{|s(n,2)|}{n!} \right).
\end{equation}

\section{Denominators with complex roots} \label{sec-catalan} 
\setcounter{equation}{0}

In this section we consider the simplest example of the type (\ref{gen-type}),
where the rational function $R_{1}(x)$ has a complex (non-real) pole. 
This is
\begin{equation}
G := - \int_{0}^{1} \frac{\ln x}{1+x^{2}} \, dx.
\nonumber
\end{equation}
\noindent
The constant $G$ is called {\em Catalan's constant} and is given by
\begin{equation}
G = \sum_{k=0}^{\infty} \frac{(-1)^{k}}{(2k+1)^{2}}.
\end{equation}
\noindent
Entry $\mathbf{4.231.12}$ of \cite{gr} states 
\begin{equation}
\int_{0}^{1} \frac{\ln x}{1+x^{2}} \, dx = -G.
\label{catalan-1}
\end{equation}

To verify (\ref{catalan-1}) simply expand the integrand to produce
\begin{eqnarray}
\int_{0}^{1} \frac{\ln x \, dx}{1+x^{2}} & =  &
- \sum_{k=0}^{\infty} (-1)^{k} \int_{0}^{1} x^{2k} \ln x \, dx \nonumber \\
 & = & -\sum_{k=0}^{\infty} (-1)^{k} \ift te^{-(2k+1)t} \, dt \nonumber \\
 & = & - \sum_{k=0}^{\infty} \frac{(-1)^{k}}{(2k+1)^{2}} \ift \sigma 
e^{-\sigma} \, d \sigma
\nonumber 
\end{eqnarray}
\noindent
The integral is evaluated by integration by parts or recognizing its value 
as $\Gamma(2)=1$. \\ 

Integration by parts gives the alternative form 
\begin{equation}
\int_{0}^{1} \frac{\tan^{-1} x}{x} \, dx  = G,
\end{equation}
\noindent
that appears as $\mathbf{4.531.1}$. \\

There are many definite integrals in \cite{gr} that are related to 
Catalan's constant. For example, the change of variables $t = 1/x$ yields from 
(\ref{catalan-1}), the value
\begin{equation}
\int_{1}^{\infty} \frac{\ln t \, dt}{1+t^{2}} = G. 
\label{formula43}
\end{equation}
\noindent
This is the second part of $4.231.12$.  Adding these two expressions for $G$,
we conclude that
\begin{equation}
\int_{0}^{\infty} \frac{\ln x \, dx}{1+x^{2}} = 0. 
\label{pole-i1}
\end{equation}

The change of variables $x = at$ in (\ref{catalan-1}) yields
$\mathbf{4.231.11}$:
\noindent
\begin{equation}
\int_{0}^{a} \frac{\ln \, dx}{x^{2}+a^{2}} = \frac{\pi \, \ln a}{4a} - 
\frac{G}{a}
\end{equation}

We now introduce material that will provide a genaralization 
of (\ref{pole-i1}) to the case of 
a multiple pole at $i$. The integral is expressed in terms of the 
{\em polygamma function}
\begin{equation}
\psi(x) = \frac{\Gamma'(x)}{\Gamma(a)}.
\end{equation}

\begin{Lem}
Let $a, \, b  \in \mathbb{R}$. Then
\begin{equation}
\int_{0}^{\infty} \frac{\ln t  \, dt }{(1+t^{2})^{b}} =
\frac{\Gamma(\tfrac{1}{2}) \psi(\tfrac{1}{2}) - 
\Gamma(b - \tfrac{1}{2}) \psi(b - \tfrac{1}{2}) }{2 \Gamma(b)}. 
\label{formula99}
\end{equation}
\end{Lem}
\begin{proof}
Define 
\begin{equation}
f(a,b) := \int_{0}^{\infty} \frac{t^{a} \, dt}{(1+t^{2})^{b}}.
\end{equation}
\noindent
Then
\begin{equation}
\frac{d}{da} f(a,b) = \int_{0}^{\infty} \frac{t^{a} \, \ln t }{(1+t^{2})^{b}}
\, dt.
\end{equation}
\noindent
In particular,
\begin{equation}
\frac{d}{da} f(a,b) \Big{|}_{a=0}  
= \int_{0}^{\infty} \frac{\ln t  \, dt }{(1+t^{2})^{b}}.
\end{equation}
\noindent
The change of variables $s = t^{2}$ gives 
\begin{equation}
f(a,b) = \frac{1}{2} \int_{0}^{\infty} \frac{s^{(a-1)/2} \, ds }{(1+s)^{b}}.
\end{equation}
\noindent
Now use the integral representation
\begin{equation}
B(u,v) = \int_{0}^{\infty} \frac{s^{u-1} \, ds}{(1+s)^{u+v}}
\end{equation}
\noindent
(given in $\mathbf{8.380.3}$ in \cite{gr} and proved in \cite{moll-gr6}) with
$u = (a+1)/2$ and $v = b - (a+1)/2$ to obtain
\begin{equation}
f(a,b) = B \left( \frac{a+1}{2}, b - \frac{a+1}{2} \right) = 
\frac{\Gamma((a+1)/2) \, \Gamma(b - (a+1)/2) }{\Gamma(b)}.
\end{equation}
\noindent
Therefore
\begin{eqnarray}
\int_{0}^{\infty} \frac{\ln t  \, dt }{(1+t^{2})^{b}} & = & 
\frac{d}{da} f(a,b) \Big{|}_{a=0}  \\ 
& = & \frac{1}{2 \Gamma(b)} \left( \Gamma'((a+1)/2) - 
\Gamma'( b - (a+1)/2) \right) \Big{|}_{a=0}. \nonumber
\end{eqnarray}
\noindent
Now use the relation $\Gamma'(x) = \psi(x) \Gamma(x)$ to obtain the result.
\end{proof}

The case $b = n \in \mathbb{N}$ requieres the value
\begin{equation}
\Gamma( n + \tfrac{1}{2}) = \frac{\sqrt{\pi}}{2^{2n}} \, 
\frac{(2n)!}{n!},
\label{gamma12}
\end{equation}
\noindent
and
\begin{equation}
\psi(n + \tfrac{1}{2}) = - \gamma + 2 \ln 2 - 2 \sum_{k=1}^{n} \frac{1}{2k-1}
= -\gamma - 2 \ln 2 + 2H_{2n} - H_{n},
\label{psi12}
\end{equation}
\noindent
that appears in $\mathbf{8.366.3}$. Here $H_{n}$ is the $n$-th harmonic 
number. The reader will find a proof of this evaluation in 
\cite{irrbook}, page 212. A proof of (\ref{gamma12}) appears as Exercise
$10.1.3$ on page 190 of \cite{irrbook}. \\

\begin{Cor}
\label{coro11}
Let $n \in \mathbb{N}$. Then
\begin{equation}
\int_{0}^{\infty} \frac{\ln x \, dx}{(1+x^{2})^{n+1}} = 
- \frac{\pi}{2^{2n+1}} \binom{2n}{n} \sum_{k=1}^{n} \frac{1}{2k-1}.
\end{equation}
\end{Cor}

We now provide a proof of Entry ${\mathbf{4.231.7}}$ in \cite{gr}:
\begin{equation}
\ift \frac{\ln x \, dx }{(a^{2}+b^{2}x^{2})^{n}} = 
\frac{\Gamma(n- \tfrac{1}{2}) \sqrt{\pi}}{4(n-1)! \, a^{2n-1}b} 
\left( 2 \ln \left( \frac{a}{2b} \right) - \gamma - \psi(n - \tfrac{1}{2}) 
\right).
\nonumber
\end{equation}

Taking the factor $b^{2}$ out of the parenthesis on the left and letting 
$c = a/b$ yields the equivalent form 
\begin{equation}
\ift \frac{\ln x \, dx }{(c^{2}+x^{2})^{n}} = 
\frac{\Gamma(n- \tfrac{1}{2}) \sqrt{\pi}}{4(n-1)!} 
\left( 2 \ln \left( \frac{c}{2} \right) - \gamma - \psi(n - \tfrac{1}{2}) 
\right).
\nonumber
\end{equation}
\noindent
It is more convenient to replace $n$ by $n+1$ to obtain
\begin{equation}
\ift \frac{\ln x \, dx }{(c^{2}+x^{2})^{n+1}} = 
\frac{\Gamma(n+ \tfrac{1}{2}) \sqrt{\pi}}{4n!} 
\left( 2 \ln \left( \frac{c}{2} \right) - \gamma - \psi(n + \tfrac{1}{2}) 
\right).
\nonumber
\end{equation}
\noindent
Using (\ref{gamma12}) and (\ref{psi12}) the requested evaluation amounts to
\begin{equation}
\ift \frac{\ln x \, dx }{(c^{2}+x^{2})^{n+1}} = 
\frac{\pi}{(2c)^{2n+1}} \binom{2n}{n} \left( \ln c - \sum_{k=1}^{n} 
\frac{1}{2k-1} \right). 
\label{int-42317}
\end{equation}
\noindent
This can be written as 
\begin{equation}
\ift \frac{\ln x \, dx }{(c^{2}+x^{2})^{n+1}} = 
\frac{\pi}{(2c)^{2n+1}} \binom{2n}{n} \left( \ln c - H_{n}+2H_{2n} \right).
\label{int-42317a}
\end{equation}
\noindent
To establish this, make the change of variables $x = ct$ to obtain
\begin{equation}
\ift \frac{\ln x \, dx }{(x^{2}+c^{2})^{n+1}} = 
\frac{\ln c}{c^{2n+1}} \ift \frac{dx }{(t^{2}+1)^{n+1}} +
\frac{1}{c^{2n+1}} \ift \frac{\ln t \, dt }{(t^{2}+1)^{n+1}}. \nonumber
\end{equation}

Using Wallis' formula
\begin{equation}
\ift \frac{dt}{(1+t^{2})^{n+1}} = \frac{\pi}{2^{2n+1}} \binom{2n}{n},
\end{equation}
\noindent
the required evaluation now follows from Corollary \ref{coro11}. 

The special case $n=0$ yields $\mathbf{4.231.8}$:
\begin{equation}
\ift \frac{\ln x \, dx }{a^{2}+b^{2}x^{2}} = \frac{\pi}{2ab} \ln \left( 
\frac{a}{b} \right). 
\end{equation}
\noindent
This evaluation also appears as $\mathbf{4.231.9}$ in the form
\begin{equation}
\ift \frac{\ln px \, dx }{q^{2}+x^{2}} = \frac{\pi}{2q} \ln pq. 
\end{equation}

\section{The case of a single purely imaginary pole} \label{sec-imag} 
\setcounter{equation}{0}

In this section we evaluate the integral
\begin{equation}
h_{n,2}(a,b) := \int_{0}^{b} \frac{\ln t \, dt}{(t^{2}+a^{2})^{n+1}},
\end{equation}
\noindent
for $a, \, b > 0$ and $n \in \mathbb{N}$. This  is 
the generalization of (\ref{h1-def}) to the case in which the 
integrand has a multiple pole at $t = ia$. The change of variables $t = ax$
yields
\begin{equation}
h_{n,2}(a,b)  = a^{-2n-1} g_{n}(b/a),
\end{equation}
\noindent
where 
\begin{equation}
g_{n}(x) := \int_{0}^{x} \frac{\ln t \, dt}{(t^{2}+1)^{n+1}}.
\end{equation}

We produce first a recurrence for the companion integral
\begin{equation}
f_{n}(x) := \int_{0}^{x} \frac{dt}{(t^{2}+1)^{n+1}}.
\label{def-f}
\end{equation}

\begin{Lem}
Let $n \in \mathbb{N}$ and $x >0$. Then
\begin{equation}
2n f_{n}(x) = (2n-1)f_{n-1}(x) + \frac{x}{(x^{2}+1)^{n}},
\label{recur-f}
\end{equation}
\noindent
and $f_{0}(x) = \tan^{-1}x$. 
\end{Lem}
\begin{proof}
Integrate by parts.
\end{proof}

The recurrence (\ref{recur-f}) is now solved using the following result 
established in \cite{moll-gr5}. 

\begin{Lem}
\label{sol-rec}
Let $n \in \mathbb{N}$ and define $\lambda_{j} = 2^{2j} \binom{2j}{j}^{-1}$.
Suppose the sequence $z_{n}$ satisfy the recurrence 
$2nz_{n}-(2n-1)z_{n-1} = r_{n}$, with $r_{n}$ given. Then
\begin{equation}
z_{n} = \frac{1}{\lambda_{n}} \left( z_{0} + \sum_{k=1}^{n} 
\frac{\lambda_{k} \, r_{k} }{2k} \right).
\end{equation}
\end{Lem}

We conclude with an explicit expression for $f_{n}(x)$. 

\begin{Prop}
\label{prop62}
Let $n \in \mathbb{N}$. Then 
\begin{equation}
\int_{0}^{x} \frac{dt}{(1+t^2)^{n+1}} = \frac{\binom{2n}{n} }{2^{2n}} 
\left( \tan^{-1}x + \sum_{j=1}^{n} \frac{2^{2j}}{2j \binom{2j}{j}} 
\frac{x}{(x^{2}+1)^{j}} \right).
\end{equation}
\end{Prop}

\noindent
\begin{Note}
This expression for $f_{n}$ appears as $\mathbf{2.148.4}$ of 
\cite{gr} in the alternative form
\begin{eqnarray}
\int_{0}^{x} \frac{dt}{(1+t^{2})^{n+1}} & =  & 
\frac{x}{2n+1} \sum_{k=1}^{n} \frac{(2n+1)!!}{(2n-2k+1)!!} \, 
\frac{(n-k)!}{2^{k} n!} \frac{1}{(1+x^{2})^{n+1-k}} \label{alter1} \\
 & + & \frac{(2n-1)!!}{2^{n} n!} \tan^{-1}x. \nonumber
\end{eqnarray}
\end{Note}

\medskip

We now produce a recurrence for the integral $g_{n}(x)$. \\

\begin{Lem}
Let $n \in \mathbb{N}$. Then the integrals $g_{n}(x)$ satisfy
\begin{equation}
2ng_{n}(x) - (2n-1)g_{n-1}(x) = 2n \ln x f_{n}(x) - 
\left[ (2n-1) \ln x + 1 \right] f_{n-1}(x).
\label{recur-g}
\end{equation}
\end{Lem}
\begin{proof}
Integration by parts yields
\begin{equation}
g_{n}(x) = f_{n}(x) \ln x - \int_{0}^{x} f_{n}(t) \frac{dt}{t}.
\label{rec111}
\end{equation}
\noindent
From the recurrence (\ref{recur-f}) we obtain
\begin{equation}
2n \int_{0}^{x} f_{n}(t) \frac{dt}{t} - 
(2n-1) \int_{0}^{x} f_{n-1}(t) \frac{dt}{t} = f_{n-1}(x). 
\end{equation}
\noindent
Now replace the integral term from (\ref{rec111}) to obtain the result.
\end{proof}

In order to produce a closed-form expression for $g_{n}(x)$ using Lemma 
\ref{sol-rec}, we need the initial condition 
\begin{equation}
g_{0}(x) = \int_{0}^{x} \frac{\ln t \, dt}{1+t^{2}}.
\end{equation}

\begin{Lem}
The function $g_{n}(x)$ is given by 
\begin{equation}
g_{n}(x) = \sum_{k=0}^{\infty} (-1)^{k} \binom{n+k}{k} 
\frac{x^{2k+1}}{2k+1} \left( \ln x - \frac{1}{2k+1} \right). 
\end{equation}
\end{Lem}
\begin{proof}
The binomial theorem gives 
\begin{equation}
(1+t^{2})^{-n-1} = \sum_{k=0}^{\infty} (-1)^{k} \binom{n+k}{k} t^{2k}.
\end{equation}
\noindent
The expression for $g_{n}(x)$ now follows by integrating term by term and
the evaluation
\begin{equation}
\int_{0}^{x} t^{2k} \ln t \, dt = \left( \ln x - \frac{1}{2k+1} \right) 
\frac{x^{2k+1}}{2k+1}.
\end{equation}
\end{proof}

\medskip

In particular, the initial condition $g_{0}(x)$ of the recurrence 
(\ref{recur-g}) is given by
\begin{equation}
g_{0}(x) = \ln x \tan^{-1}x - L(x),
\end{equation}
\noindent
where
\begin{equation}
L(x) = \sum_{k=0}^{\infty} \frac{(-1)^{k} x^{2k+1}}{(2k+1)^{2}} =
\int_{0}^{x} \frac{\tan^{-1}t}{t} \, dt. 
\end{equation}

\medskip

The recurrence (\ref{recur-g}) is now solved using (\ref{sol-rec}) to 
produce
\begin{eqnarray}
g_{n}(x) & = & 2^{-2n} \binom{2n}{n} g_{0}(x)  \nonumber  \\
 & + & 
\frac{\binom{2n}{n}}{2^{2n}} 
\sum_{j=1}^{n} \frac{2^{2j}}{2j \binom{2j}{j} } 
\left[ \left\{  2j f_{j}(x) - (2j-1)f_{j-1}(x) \right\} \ln x - f_{j-1}(x) \right]. 
\nonumber
\end{eqnarray}
\noindent
Using the recurrence (\ref{recur-f}), this can be written as 
\begin{equation}
g_{n}(x) = \frac{\binom{2n}{n}}{2^{2n}} 
\left( g_{0}(x) + 
\sum_{j=1}^{n} \frac{2^{2j}}{2j \binom{2j}{j} } 
\left[ \frac{x \, \ln x }{(x^{2}+1)^{j}} - f_{j-1}(x) \right] \right). 
\end{equation}

Solving the recurrence yields: \\

\begin{Thm}
\label{thm-g}
Let $n \in \mathbb{N}$. Define the rational function
\begin{equation}
p_{j}(x) = 
\sum_{k=1}^{j} \frac{2^{2k}}{2k \binom{2k}{k}} \frac{x}{(1+x^{2})^{k}}
\end{equation}
\noindent
Then the integral $g_{n}(x)$ is given by
\begin{equation}
\int_{0}^{x} \frac{\ln t \, dt}{(1+t^{2})^{n+1}}  
= \frac{\binom{2n}{n}}{2^{2n}} \left[ g_{0}(x) + p_{n}(x) \ln x - 
\sum_{k=1}^{n} 
\frac{\tan^{-1} x  + p_{k-1}(x)}{2k-1} \right]. 
\nonumber
\end{equation}
\end{Thm}

\medskip

\noindent
\begin{Note}
The special case $x=1$ in (\ref{thm-g}) produces 
\begin{equation}
\int_{0}^{1} \frac{\ln t \, dt}{(t^{2}+1)^{n+1}} = 2^{-2n} \binom{2n}{n} 
\left( G - \sum_{k=1}^{n} \frac{ \tfrac{\pi}{4} + p_{k-1}(1)}{2k-1} \right).
\end{equation}
The values 
\begin{equation}
p_{n}(1) = \frac{1}{2} \sum_{j=1}^{n} \frac{2^{j}}{j \binom{2j}{j} }
\end{equation}
\noindent
do not admit a closed-form, but they do satisfy the three term recurrence 
\begin{equation}
(2n+1)p_{n+1}(1) - (3n+1)p_{n}(1) + n p_{n-1}(1) = 0. 
\end{equation}
\noindent
The reader is invited to verify the expansion
\begin{equation}
\sum_{k=1}^{\infty} \frac{x^{k}}{k \binom{2k}{k}} = 
\frac{2 \sqrt{x} \, \sin^{-1}(\sqrt{x}/2) }{\sqrt{4-x}},
\end{equation}
\noindent
from which it follows that
\begin{equation}
\sum_{k=1}^{\infty} \frac{2^{k}}{k \, \binom{2k}{k}} = \frac{\pi}{2}.
\end{equation}
\end{Note}

\medskip

\noindent
{\bf An alternative derivation}. Integration by parts produces 
\begin{equation}
\int_{0}^{a} \frac{\ln s \, ds}{s^{2}+b} = 
\frac{1}{\sqrt{b}} \ln a \tan^{-1} \frac{a}{\sqrt{b}} - 
\frac{1}{\sqrt{b}} \int_{0}^{a/\sqrt{b}} \frac{\tan^{-1}x}{x} \, dx. 
\end{equation}

Differentiating $n$-times with respect to $b$ and using 
\begin{equation}
\left( \frac{d}{db} \right)^{j} \frac{1}{\sqrt{b}}  = 
\frac{(-1)^{j} \, (2j)!}{j! 2^{2j} b^{j+1/2}}, \,
\left( \frac{d}{db} \right)^{j} \frac{1}{a^{2}+b}  = 
\frac{(-1)^{j} \, j!}{(a^{2}+b)^{j+1}}, \nonumber
\end{equation}
\noindent
and 
\begin{equation}
\left( \frac{d}{db} \right)^{j} \tan^{-1} \frac{a}{\sqrt{b}} = 
(-1)^{j} (j-1)! a 
\sum_{k=0}^{j-1} \frac{\binom{2k}{k}}{2^{2k+1} b^{k+1/2} (a^{2}+b)^{j-k}},
\nonumber
\end{equation}
we obtain
\begin{eqnarray}
\int_{0}^{a} \frac{\ln s \, ds}{(s^{2}+b)^{n+1}} & = & 
\ln a \tan^{-1} \left( \frac{a}{\sqrt{b}} \right) F_{n}(b) + 
\frac{a \ln a }{2} \sum_{k=1}^{n} \frac{F_{n-k}(b)}{k} 
\sum_{j=0}^{k-1} \frac{F_{j}(b)}{(a^{2}+b)^{k-j}} \nonumber \\
& - & F_{n}(b) \int_{0}^{a} \frac{\tan^{-1}x}{x} dx - 
\frac{1}{2} \sum_{k=1}^{n} \frac{F_{n-k}(b)}{k} 
\sum_{j=0}^{j-1} F_{j}(b) \int_{0}^{a} \frac{dt}{(t^{2}+b)^{k-j}} 
\nonumber
\end{eqnarray}
\noindent
with $F_{j}(b) = 2^{-2j}b^{-j-1/2}\binom{2j}{j}$. The last integral in this 
expression can be evaluated using Proposition \ref{prop62} to produce
a generalization of Theorem \ref{thm-g}. We have replaced the parameter $b$
by $b^{2}$ to produce a cleaner formula. \\

\begin{Thm}
Let $a, \, b \in \mathbb{R}$ with $a>0$ and $n \in \mathbb{N}$. Introduce the 
notation
\begin{equation}
 F_{j}(b) = \frac{\binom{2j}{j}}{2^{2j} b^{2j+1}}. 
\end{equation}
\noindent
Then 
\begin{eqnarray}
\int_{0}^{a} \frac{\ln s \, ds}{(s^{2}+b^{2})^{n+1}} & = & 
F_{n}(b) \ln a \tan^{-1}(a/b) + \frac{a \ln a }{2} 
\sum_{k=1}^{n} \frac{F_{n-k}(b)}{k} \sum_{j=0}^{k-1} \frac{F_{j}(b)}
{(a^{2}+b^{2})^{k-1}} \nonumber \\
& - & F_{n}(b) \int_{0}^{a/b} \frac{\tan^{-1}x}{x} \, dx - 
\frac{1}{2}\tan^{-1}(a/b) \left( 
\sum_{k=1}^{n} \frac{F_{n-k}(b)}{k} \sum_{j=0}^{k-1} F_{j}(b) F_{k-j-1}(b) 
\right) \nonumber \\
& - & \frac{a}{4 \sqrt{b}} \sum_{k=1}^{n} \frac{F_{n-k}(b)}{k} 
\sum_{j=0}^{k-1} F_{j}(b) F_{k-j-1}(b) 
\sum_{r=1}^{k-j-1} \frac{1}{rF_{r}(b) (a^{2}+b^{2})^{r}}. \nonumber
\end{eqnarray}
\end{Thm}

\section{Some trigonometric versions} \label{sec-trigo} 
\setcounter{equation}{0}

In this section we provide trigonometric versions of some of the evaluations
provided in the previous sections. Many of these integrals correspond to 
special values of the {\em Lobachevsky function} defined by
\begin{equation}
L(x) := - \int_{0}^{x} \ln \cos t \, dt.
\label{lcos}
\end{equation}
\noindent
This appears as entry $\mathbf{8.260}$ in \cite{gr} and also as 
$\mathbf{4.224.4}$. The change of variables $t = \pi/2 - x$ provides a 
proof of $\mathbf{4.224.1}$:
\begin{equation}
\int_{0}^{x} \ln \sin t \, dt = L \left(\pi/2 -x \right) - 
L \left( \pi/2 \right). 
\label{lsin}
\end{equation}
\noindent 
The special value 
\begin{equation}
L \left( \pi/2 \right) = \int_{0}^{\pi/2} \ln \sin x \, dx = 
- \frac{\pi}{2} \ln 2,
\end{equation}
\noindent
appears as $\mathbf{4.224.3}$.  The 
change of variables $t = \tfrac{\pi}{2} - x$ yields $\mathbf{4.224.6}$:
\begin{equation}
\int_{0}^{\pi/2} \ln \cos x \, dx = 
- \frac{\pi}{2} \ln 2.
\label{formula54}
\end{equation}

To establish these evaluations, observe that, by 
symmetry,
\begin{eqnarray}
2 \int_{0}^{\pi/2} \ln \sin x \, dx & = & 
\int_{0}^{\pi/2} \ln \sin  x \, dx + 
\int_{0}^{\pi/2} \ln \cos  x \, dx \nonumber \\
 & = & \int_{0}^{\pi/2} \ln \left( \frac{\sin  2x}{2} \right) \, dx 
\nonumber \\
 & = & \int_{0}^{\pi/2} \ln(\sin  2x) \, dx  - \frac{\pi}{2} \ln 2.
\nonumber
\end{eqnarray}
\noindent
Now let $t = 2x$ in the last integral to obtain the result.  \\

Combining (\ref{lcos}) and (\ref{lsin}) we obtain $\mathbf{4.227.1}$
\begin{equation}
\int_{0}^{u} \ln \tan x \, dx  = L(u) + L(\pi/2-u) + \frac{\pi}{2} \ln 2.
\end{equation}

\medskip

The identity (\ref{trig1}) and the evaluation (\ref{catalan-1}) yield
the value of $\mathbf{4.227.2}$:
\begin{equation}
\int_{0}^{\pi/4} \ln \tan t \, dt = -G.
\end{equation}
\noindent
Now observe that
\begin{equation}
\int_{0}^{\pi/4} \ln \tan t \, dt = 
\int_{0}^{\pi/4} \ln \sin t \, dt  - 
\int_{0}^{\pi/4} \ln \cos t \, dt  = -G 
\end{equation}
\noindent
and 
\begin{equation}
\int_{0}^{\pi/2} \ln \sin t \, dt = 
\int_{0}^{\pi/4} \ln \sin t \, dt  + 
\int_{0}^{\pi/4} \ln \cos t \, dt  = - \frac{\pi}{2} \ln 2. 
\end{equation}
\noindent
Solving this system of equations yields
\begin{equation}
\int_{0}^{\pi/4} \ln \sin t \, dt = - \frac{\pi}{4} \ln 2 - \frac{G}{2}
\label{formula59}
\end{equation}
\noindent 
that appears as $\mathbf{4.224.2}$ in \cite{gr} and 
\begin{equation}
\int_{0}^{\pi/4} \ln \cos t \, dt = - \frac{\pi}{4} \ln 2 + \frac{G}{2}
\label{formula510}
\end{equation}
\noindent 
that appears as $\mathbf{4.224.5}$.  \\

We now make use of the identity 
\begin{equation}
\cos x - \sin x = \sqrt{2} \cos(x + \pi/4)
\end{equation}
\noindent
to obtain
\begin{eqnarray}
\int_{0}^{\pi/4} \ln( \cos x - \sin x) \, dx & =  & 
\frac{\pi}{8} + \int_{0}^{\pi/2} \ln \cos( x + \pi/4) \nonumber \\
& = & \frac{\pi}{8} \ln 2 + \int_{0}^{\pi/2} \ln \cos t \, dt - 
\int_{0}^{\pi/4} \ln \cos t \, dt. \nonumber
\end{eqnarray}
\noindent
The first integral is given in (\ref{formula54}) as $- \frac{\pi}{2} \ln 2$
and the second one as $- \frac{\pi}{4} \ln 2 + \frac{G}{2}$ in 
(\ref{formula510}). We conclude with the evaluation of $\mathbf{4.225.1}$
\begin{equation}
\int_{0}^{\pi/4} \ln( \cos x - \sin x ) \, dx = 
-\frac{\pi}{8} \ln 2 - \frac{G}{2}.
\end{equation}
\noindent
A similar analysis produces $\mathbf{4.225.2}$:
\begin{equation}
\int_{0}^{\pi/4} \ln( \cos x + \sin x ) \, dx = 
-\frac{\pi}{8} \ln 2 + \frac{G}{2}.
\label{4.225.2}
\end{equation}
These evaluations can be used to produce $\mathbf{4.227.9}$:
\begin{equation}
\int_{0}^{\pi/4} \ln( 1 + \tan  x ) \, dx = 
\frac{\pi}{8} \ln 2, 
\end{equation}
\noindent
and $\mathbf{4.227.11}$:
\begin{equation}
\int_{0}^{\pi/4} \ln( 1 - \tan  x ) \, dx = 
\frac{\pi}{8} \ln 2 - G. 
\end{equation}
\noindent
To prove these formulas, simply write
\begin{equation}
\ln(1 \pm  \tan x ) = \ln(\cos \pm \sin x) - \ln \cos x.
\end{equation}
\noindent
The same type of calculations provide verification of $\mathbf{4.227.13}$
\begin{equation}
\int_{0}^{\pi/4} \ln( 1 + \cot x ) \, dx = \frac{\pi}{8} \ln 2 + G 
\end{equation}
\noindent
and $\mathbf{4.227.14}$
\begin{equation}
\int_{0}^{\pi/4} \ln( \cot x -1 ) \, dx = \frac{\pi}{8} \ln 2. 
\end{equation}
\noindent
The next example of this type is $\mathbf{4.227.15}$:
\begin{equation}
\int_{0}^{\pi/4} \ln( \tan x + \cot x ) \, dx = \frac{\pi}{2} \ln 2, 
\end{equation}
\noindent
This is evaluated by writing the integral as 
\begin{equation}
- \int_{0}^{\pi/4} \ln( \sin  x) \, dx - \int_{0}^{\pi/2} 
\ln( \cos x) \, dx = \frac{\pi}{2} \ln 2, 
\end{equation}
\noindent
using (\ref{formula59}) and (\ref{formula510}). 

The evaluation of $\mathbf{4.227.10}$ 
\begin{equation}
\int_{0}^{\pi/2} \ln( 1 + \tan x ) \, dx = \frac{\pi}{4} \ln 2 + G,
\end{equation}
\noindent
follows from the integrals evaluated here. Indeed,
\begin{eqnarray}
\int_{0}^{\pi/2} \ln( 1 + \tan x ) \, dx & =  & 
\int_{0}^{\pi/2} \ln( \sin x  + \cos x ) \, dx - 
\int_{0}^{\pi/2} \ln(  \cos x ) \, dx  \nonumber \\
& = & 2 \int_{0}^{\pi/4} \ln( \sin x  + \cos x ) \, dx - 
\int_{0}^{\pi/2} \ln(  \cos x ) \, dx  \nonumber \\
& = & 2 \left( - \frac{\pi}{8} \ln 2 + \frac{G}{2} \right) 
+ \frac{\pi}{2} \ln 2, 
\nonumber
\end{eqnarray}
\noindent
where we have used (\ref{formula54}) and (\ref{4.225.2}).

\medskip

The identity (\ref{pole-i1}) yields 
\begin{equation}
\int_{0}^{\pi/2} \ln \tan t \, dt = 0. 
\end{equation}
\noindent
The apparent generalization 
\begin{equation}
\int_{0}^{\pi/2} \ln( a  \tan t) \, dt = \frac{\pi}{2} \ln a,
\end{equation}
\noindent
with $a> 0$, appears as $\mathbf{4.227.3}$.

\medskip

The evaluations (\ref{formula59}) and (\ref{formula510}) can be brought 
back into rational form. The change of variables $t = \tan^{-1}u$ produces 
from (\ref{formula510}):
\begin{equation}
-\frac{\pi}{4} \ln 2 + \frac{G}{2}  =  
- \frac{1}{2} \int_{0}^{1} \frac{\ln(1+u^{2})}{1+u^{2}} \, du. 
\end{equation}
\noindent
We have obtained a proof of $\mathbf{4.295.5}$:
\begin{equation}
\int_{0}^{1} \frac{\ln(1+x^{2})}{1+x^{2}} \, dx = \frac{\pi}{2} \ln 2 - G.
\end{equation}
\noindent
The change of variables $t = 1/x$ and (\ref{formula43}) yield 
$\mathbf{4.295.6}$:
\begin{equation}
\int_{1}^{\infty} \frac{\ln(1+t^{2})}{1+t^{2}} \, dt = \frac{\pi}{2} \ln 2 
+ G.
\end{equation}
\noindent

\medskip 

There are many other intergrals that may be evaluated by the methods reported
here. For instance, integration by parts yields 
\begin{equation}
\int_{0}^{x} t \, \cot t \, dt  = x \, \ln \sin x - 
\int_{0}^{} \ln \sin t \, dt. 
\end{equation}
\noindent
Using (\ref{lsin}) we obtain 
\begin{equation}
\int_{0}^{x} t \, \cot t \, dt = x \, \ln \sin x - 
L \left( \frac{\pi}{2} - x \right)
+ \frac{\pi}{2} \ln 2. 
\end{equation}
\noindent
In particular, from $L(0)=0$, we obtain $\mathbf{3.747.7}$:
\begin{equation}
\int_{0}^{\pi/2} \frac{t \, dt}{\tan t } = \frac{\pi}{2} \ln 2. 
\end{equation}
\noindent 
The change of  variables $u = \sin x$ produces from here the evaluation of 
$\mathbf{4.521.1}$:
\begin{equation}
\int_{0}^{1} \frac{\text{Arcsin }u}{u} \, du = \frac{\pi}{2} \ln 2. 
\end{equation}

Further evaluations will be reported in a future publication.

\end{document}